\input amstex
\documentstyle{amsppt}

\magnification=\magstep1
\def\bb#1{{\bold #1}}
\def\pp#1{\bb P^{#1}}
\def\pas{\par\noindent}
\def\expdim{\operatorname{expdim}}

\topmatter
\author L.Chiantini, M.Coppens \endauthor
\title Grassmannians of secant varieties \endtitle
\abstract
For an irreducible projective variety $X$, we study the family of $h$-planes contained in the secant variety Sec${}_k(X)$, for $0<h<k$. These families have an expected dimension and we study varieties for which the expected dimension is not attained; for these varieties, making general consecutive projections to lower dimensional spaces, we do not get the expected singularities. In particular, we examine the family $G_{1,2}$ of lines sitting in $3$-secant planes to a surface $S$. We show that the actual dimension of $G_{1,2}$ is equal to the expected dimension unless $S$ is a cone or a rational normal scroll of degree $4$ in $\pp 5$.
\endabstract
\endtopmatter

\heading 0. Introduction\endheading
Let $X\subset\pp r$ be a smooth, non degenerate $n$-dimensional projective variety. Take $P\in\pp r$ and assume that for all lines $L$ through $P$ the intersection, as a scheme, has length at most one. Then using the projection with center $P$ one obtains an embedding $X\subset\pp {r-1}$. Such a point $P$ exists if and only if $S_1(X)\neq\pp r$, with $S_1(X)$ being the secant variety of $X$: it is the closure of lines spanned by pairs of distinct points in $X$ (see [JH] lecture 15). Clearly $\dim(S_1(X))\le 2n+1$ and $\min\{r,2n+1\}$ can be considered as the ''expected dimension'' of $S_1(X)$. Hence, using such projections, we can embed $X$ in $\pp {2n+1}$; we can go further and project $X$ isomorphically in some $\pp  m$, $m<2n+1$ if and only if $S_1(X)$ has dimension smaller than the expected one.\smallskip

The classification of varieties for which $S_1(X)$ has dimension less than the expected value was studied by Severi, Terracini and Scorza (see e.g. [Sc]) for objects of small dimension and it has been recently reconsidered by several authors. Severi found that the Veronese surface is the unique smooth surface in $\pp r$, $r\ge 5$ that can be projected isomorphically to $\pp 4$. In [Z], lower bounds for $\dim(S_1(X))$ are proved and a classification for varieties attained this lower bound is presented. These varieties can be projected isomorphically to some projective space of dimension much smaller than $2n+1$.\par
\smallskip

In general, projecting a variety $X\subset\pp {2n+1}$ from some disjoint linear subspace $\pi$ of dimension $k>0$, as $k$ increases  one expects that points of higher multiplicity must arise in the image. A way (unfortunately not the unique one) to obtain these  multiple points is by considering a linear span $\pi'$ of $k+2$ {\it distinct points} of $X$ which contains $\pi$: clearly $\pi'$ is contracted to a $(k+2)-uple$ point of the projection of $X$. The existence of $\pi'$ for a general choice of the center $\pi$ of projection can be rephrased as follows: consider the Grassmannian $G(k,r)$ of $k$-planes in $\pp r$ and for a general set $P_1,\dots,P_{k+2}$ of points of $X$ consider the subset of $G(k,r)$ formed by $k$-planes contained in the span of the point $P_i$'s. As the $P_i$'s move, these subsets describe an algebraic subspace $G(X)\subset G(k,r)$ which has an obvious expected dimension (see section 1 for more precise definitions) and the singularity above occurs when the closure of $G(X)$ coincides  with $G(k,r)$.\pas
In this setting, varieties for which some space $G(X)$ above has dimension less than the expected one are the analogous of varieties with degenerate secant variety; thus we may expect that they are strongly characterized by this property and study their classification. 
\pas Notice that this problem was partly considered by classical geometers as a generalization of the Waring problem for forms. Indeed, for varieties $X$ which are image of projective spaces under some Veronese embeddings, it arises naturally when one tries to write a set of forms  as sum of powers of the same linear forms (see [B] and [T] for wider discussion).\smallskip

In the present paper, we propose a systematic study of the subspaces of the Grassmannians arising as above  from secant spaces to a projective variety $X$; we call these spaces {\it Grassmannians of secant varieties}. In the first section we give precise (and more general) definitions. We prove then some general results on the dimensions of these spaces and we show that an exceptional behaviour gives rise to exceptional behaviour also with respect to $S_k(X)$, the natural generalization of $S_1(X)$.\smallskip 

In the second section, we obtain the classification of irreducible surfaces $X\subset\pp r$, $r\ge 5$, for which the Grassmannian of lines contained in 3-secant planes has dimension smaller than the expected one.\par
Indeed, for such surfaces one expects that lines contained in some 3-secant plane describe a subvariety of dimension 8 in $G(1,r)$; e.g. when $r=5$ one expects that a general line lies in some 3-secant plane (or, in other words, the projection of $X$ to $\pp 3$ should have some triple point). On the other hand it is classically known that the embedding of $\pp 1\times\pp 1$ in $\pp 5$ by a divisor of type $(2,1)$ does not enjoy this property. We are able to prove that this is the unique smooth surface for which the Grassmannian of lines in 3-secant planes has dimension less than the expected. Also taking singular (but irreducible) surfaces into account, we have no further examples of degenerate Grasmmannian of 3-secant planes, except for cones.\smallskip

Notice the curious behaviour of smooth surfaces of minimal degree in $\pp 5$. There are 2 types of them: scrolls and the Veronese surface. The second ones have a degenerate secant variety $S_1(X)$, but all the Grassmannians of secant varieties are as expected. Surfaces of the first type, instead, have a nice behaviour with respect to secant varieties, but their Grassmannian of 3-secant planes degenerates.
\smallskip

The authors are members of the european AGE project. The first author is supported by the italian MURST fund; he would like to thank Ciro Ciliberto, for many fruitful discussions, and the University of Leuwen, for its warm hospitality. The second author is supported by a Research Fellowship at the University of Leuwen; he would like to thank the University of Siena for its hospitality during the preparation of the manuscript.

\heading 1. General Properties \endheading

{\bf Notation}\smallskip
\pas
We work over the complex field $\bb C$.\pas
Let $X\subset\pp r$ be an integral non-degenerate variety of dimension $n$.\pas
For $k\le r$, a general $(k+1)$-uple of points in $X$ spans a $k$-plane. Therefore we have the ''span'' rational map:
$$\Phi: X^{k+1} \cdots\to G(k,r)$$
to the Grassmannian of $k$-planes in $\pp r$.\pas
For all $h<k$, consider the ''incidence'' diagram:
$$ \CD  I @>\alpha_h>> G(h,r) \\
          @V\beta_hVV            @. \\
         G(k,r)        @.
\endCD $$
where $I$ is the ''incidence relation'' of pairs $(h,H)$ with $h \subset H$.\pas
Call:\par
$G_k(X)=$ closure of the image of $\Phi$;\par
$S_k(X)=$ closure of $\alpha_0(\beta_0^{-1}(G_k(X)))$;\pas
and more generally:\par
$G_{h,k}(X)=$ closure of $\alpha_h(\beta_h^{-1}(G_k(X))) \subset G(h,r)$.\pas
In fact:\par
$G_{h,k}(X)$ = closure of $\{ h$-planes contained in some $k-$plane, $(k+1)-$secant to $X\}$; \par
$S_k(X)=G_{0,k}(X)=$  closure of the union of all  $k-$planes, $(k+1)$-secant to $X$. \smallskip
\pas
Observe that we do not consider singular points of $X$; e.g. a general line through a double point needs not to be secant in our definition.\smallskip
\pas
For the expected dimensions of these objects, we have:\par
$\expdim(G_k(X))=\min\{n(k+1), (r-k)(k+1)\};$\par
$\expdim(G_{h,k}(X))=\min\{(k-h)(h+1)+n(k+1), (r-h)(h+1)\}$ \pas
so that, as usual:\par
$\expdim(S_k(X))=\min\{n(k+1)+k, r\}.$\pas
Notice that all these varieties are irreducible, since $X$ is.

\proclaim{Proposition 1.1} $\dim G_k(X)=\expdim(G_k(X))$.
\endproclaim
\demo{Proof} First assume $n>r-k$; then any $k$-plane $\pi$ meets $X$ in (at least) a curve; moving generically $k+1$ points of this curve, we see that $\pi\in G_k(X)$. It follows that $G_k(X)$ is the whole Grassmannian, hence its dimension is $(k+1)(r-k)$.\pas
Assume now $n+k\le r$, so that $\expdim G_k(X)=n(k+1)$. Since $X^{(k+1)}$ has dimension $n(k+1)$, the actual dimension of $G_k(X)$ is always less or equal than the expected one, and equality means that $\Phi$ has finite general fibers.\pas
Assume $\dim G_k(X)<\expdim(G_k(X))$ and take $k$ minimal, with this property.  Since $k$ is minimal, the span of $k$ general points of $X$ meets $X$ in a finite set; hence the projection $X'\subset\pp {r-k+1}$ of X from $k-1$ general points still has dimension $n$, furthermore a general projection from some point of $X'$ contracts curves on $X'$; this is possible only if $X'$ is linear: a contradiction. 
\qed\enddemo

\proclaim{Corollary 1.2} $\dim G_{h,k}(X)\le \expdim G_{h,k}(X)$.\pas
$\dim G_{h,k}(X)=(k-h)(h+1)+n(k+1)$ if and only if a general $h-$plane in $G_{h,k}(X)$ lies only in a finite set of $k$-planes in $G_k(X)$.
\endproclaim
\demo{Proof} 
Immediate from the fact that $\beta_h$ has fibers of dimension $(k-h)(h+1)$.
\qed\enddemo
  
\proclaim{Corollary 1.3} If $k<r-n$ and $\pi\in G_k(X)$ is general, then $\pi\cap X$ is formed by exactly $k+1$ points.\endproclaim
\demo{Proof} Induction on $k$. If $k=1$ then this is the well known 3-secant lemma: not every secant is a 3-secant. For $k>1$, take projection from a general point of $X$.
\qed\enddemo

We recall some well-known facts.

\proclaim{Proposition 1.4}  $G_{h,k}(X)=G_{h,k+1}(X)$ implies $S_k(X)=\pp r$. \endproclaim
\demo{Proof} The condition implies $S_{k+1}(X)=S_k(X)$; then use [Z], V.1.3.
\qed\enddemo 

\proclaim{Terracini's Lemma} For general points $P_0,\dots, P_k\in X$ and $u$ general in their span, one has
$$T_{u,S_k(X)} = <T_{P_0,X}, \dots, T_{P_k,X}>.$$
\endproclaim

In fact, Terracini's lemma also works when $X$ is reducible. 

\proclaim{Linear Lemma} Any set of $m$-planes such that any two of them meet in a $(m-1)$-plane, either is contained in some fixed $\pp{m+1}$ or has a ${m-1}$-plane for base locus.
\endproclaim
\demo{Proof}  Call $H$ the $(m+1)$-plane spanned by two elements $A,B$ of the family. Assume that some $C$ in the set does not lie in $H$; then $C$ meets $H$ in  a $(m-1)$-plane,  since it meets $A,B$ in a $(m-1)$-plane, and  $C\cap H= A\cap B$; any element of the set contained in $H$ must then contain $C\cap H$; hence $A\cap B$ is the base locus.
\qed\enddemo

Using the linear Lemma and Terracini's Lemma, we can look at the situation of  secant varieties for curves.
 
\proclaim{Proposition 1.5} Let $X\subset \pp r$ be a non degenerate, reduced (but possibly reducible) curve, such that $\dim S_k(X)\le k+1<r$. Then $X$ is a cone.
\endproclaim
\demo{Proof} First take $k=1$ and assume $\dim S_1(X)<3$. Then, by Terracini's Lemma, the span of two general tangent lines to $X$ is a plane, so any pair of tangent lines meets. By [H] IV.3.8, if $X$ is non degenerate and irreducible, this is impossible. So all components of $X$ are degenerate.\pas
If $X_1$ and $X_2$ are two components, not contained in the same plane, consider a plane $\pi$ through $X_1$; a general tangent line to $X_2$ meets $\pi$ in a point $P$, which must lie in any tangent line to $X_1$; it follows that $X_1$ is strange, hence a line through $P$. Changing $X_1$ and $X_2$, we see that all the components of $X$ are lines, any two of them meeting at some points. The conclusion now follows from the Linear Lemma.\pas
For $k>1$, just work by induction: if $X$ is not a cone, then $\dim S_{k-1}(X)\ge k+1$ so $\dim S_k(X)=k+1$ implies, by Terracini's Lemma, that for a general choice of $P_1,\dots, P_k\in X$, the linear span of the tangent lines at the points $P_i$ is a $\pp{k+1}$ containing the tangent line to any other general point. This is impossible, for $X$ is non degenerate.\qed\enddemo

Next, we prove some general results on the behaviour of grassmannians of secant varieties. We are going to use them for the case of surfaces and hope they will prove useful also in higher dimensions, for the classification of varieties whose grassmannians have a degenerate behaviour.

\proclaim{Proposition 1.6} If $\dim G_{h,k}(X)=(k-h)(h+1)+n(k+1)-x$ for some $x>0$, then $\dim G_{h-1,k}(X)\le(k-h+1)h+n(k+1)-x-1$. In particular
$\dim S_k(X)\le n(k+1)+k-x-h$.
\endproclaim
\demo{Proof} The first inequality implies that $\alpha_h$ has $x$-dimensional general fiber, i.e. any $h$-plane  in some $(k+1)$-secant $k$-plane is in fact contained in a $x$-dimensional family of such planes. If $\pi'$ is a general $(h-1)$-plane in some $H\in G_k(X)$, then we have a ($k-h$)-dimensional family of $h$-planes in $H$ containing $\pi'$; the inverse image $I'$ of this family in $\alpha_h$ has dimension $k-h+x$; assume that $\dim\beta_h(I')=x'\le x$. Then over $H'\in\beta_h(I')$ general we have a fiber of dimension at least $k-h+x-x'$; this is only possible if $x=x'$ and all the $h-$planes in $H$, containing $\pi'$, are in fact contained in $H'$; it follows $H=H'$, which contradicts $x>0$. Therefore $\beta(I')$ has dimension at least $x+1$ and we are done.
\qed\enddemo

\proclaim{Proposition 1.7} Assume $S_{k-1}(X)\neq\pp r$. If $X$ is not a cone and $\dim G_{h,k}(X)=(k-h)(h+1)+n(k+1)-x$, $x>0$ then either $\dim G_{h-1,k}(X)\le  (k-h+1)h+n(k+1)-x-2$ or $\dim S_k(X)\le n(k+1)+k-k(h+x)$. In particular $\dim S_k(X)\le n(k+1)+k-x-2h$. 
\endproclaim
\demo{Proof} Fix $H\in G_k(X)$ general and let $\pi'$ be a general $(h-1)$-plane in $H$. Call $G'$ the 
$(k-h)$-dimensional family of $h$-planes in $H$ through $\pi'$ and call $I'$ its inverse image in $I$. We know that $\dim I' =k-h+x$ and its image in $G(k,r)$ is at least $(x+1)$-dimensional, with equality when $\dim G_{h-1,k}(X)=(k-h+1)h+n(k+1)-x-1$; we assume that equality holds and prove that $\dim S_k(X)\le n(k+1)+k-h(h+x)$.\pas
Take $H'$ general in the image of $I'$. The fiber of $I'$ over $H'$ is a $(k-h-1)$-dimensional family of $h$-planes $\pi$ satisfying $\pi'\subset\pi\subset H$ and $\pi\subset H'$; this implies $\dim(H\cap H')=k-1$;  $H,H'$ are general moreover the intersections of all elements of $\beta(I')$ cannot be a fixed $(k-1)$-plane, for there is an element in $\beta(I')$ through a general $h$-planes of $H$ containing $\pi'$; it follows by the Linear Lemma that the elements of $\beta(I')$ lie in a fixed $(k+1)$-plane $V$.\pas
Since $H\cap X$ generate $H$,  then $H'\cap X\neq H\cap X$ for $H'\in\beta(I')$ general; so we may write $H\cap X=T_1\cup T_2$ with $T_1\subset H'$ and $T_2\cap H'=\emptyset$ for $H'\in\beta(I')$ general. Put $t=\deg(T_2)$, so $k+1-t=\deg(T_1)$. The points of $T_2$ move when $H$ moves in $\beta(I')$; this implies that $V$ cuts $X$ in a curve $C$ of degree at least $t$. \pas
We claim that $t\ge 3$. Indeed $t>0$ and if $t<3$, then the spaces $H'$ in $\beta(I')$ have a common ($k-2$)-plane spanned by $k-1$ points of $T_1$; since $\pi'$ is general in $H$, it is not contained in this ($k-2$)-plane; since all $H'\in\beta(I')$ contain $\pi'$, this implies that they have a common $(k-1)$-plane; we yet know that this leads to a contradiction. Observe that if $C$ spans a $(t-1)$-plane, since $C\not\subset H$ then $C\cap H$ spans a $(t-2)$-plane at most and the span of $C\cap H$ and $T_1$ has dimension $\le t-2+k+1-t =k-1$; but $T_2\subset C\cap H$ and $<T_1\cup T_2>=H$, a contradiction. Thus $\dim <C>\ge t\ge 3$.\pas
We claim that $S_{t-2}(C)\neq <C>$. Otherwise the span of $S_{t-2}(C)\cup T_1$ contains $<C>\cup T_1$ hence also $ H$; but 
$<S_{t-2}(C)\cup T_1>\subset S_{t-2+k+1-t}(X)=S_{k-1}(X)$, hence $\pi'\subset S_{k-1}(X)$ which in turn implies $S_k(X)\subset S_{k-1}(X)$, i.e. $S_{k-1}(X)=\pp r$, a contradiction.\pas
Now assume $\dim S_1(C)\ge 3$; this yields $S_{\dim<C>-2}(C)=<C>$; since $H\cap C$ contains at least $\dim<C>$ points, we may find $k+1-\dim<C>$ points in $T_1$ to conclude that $V\subset S_{k-1}(X)$, hence again $S_k(X)\subset S_{k-1}(X)$, a contradiction. \pas
So we have  $\dim S_1(C)< 3$. By Proposition 1.5, we see that $C$ is a cone; we claim that it has degree $h+1+x$. Indeed the $x$-dimensional family of $k$-planes in $V$ through a general $h$-plane $\pi\supset\pi'$ must contain $k+1-\deg(C)$ fixed points in $T_1$, hence these $k$-planes contain a fixed linear subspace of $H$ of dimension $h+k+1-\deg(C)$. It follows $x=k-(h+k+1-\deg(C))$.\pas
Call $T$ the vertex of the cone $C$. Since $H$ is general, by Corollary 1.3 $H\cap X$ contains exactly $k+1$ points, so $T\notin H$; moreover, if $T$ is fixed when we move generically one point of $H\cap X$ and fix the others, then $X$ is a cone, absurd. So $T$ describes a subvariety of $\Cal T$ of $X$. We claim that $\dim \Cal T$ is a positive multiple of $(h+1+x)$; indeed consider the correspondence $Z\subset X^{h+1+x}\times X$, $Z=\{(P_1,\dots,P_{h+1+x}, T): <P_i,T>\subset X$ for all $i\}$. Then the dimension of the fiber of $Z$ over $T$  is a multiple of $x+1+h$, for all points of $X\cap H$ can be interchanged, but this fiber has also dimension $n(h+1+x)-\dim\Cal T$. In particular $\dim\Cal T \ge h+1+x$.\pas 
Fix 2 points $A,B\in H\cap X$ and let the other vary: the corresponding points $T$ describe  a subvariety of dimension $\ge h-1+x$ in $X$; this yields immediately $\dim T_A\cap T_B\ge h-1+x$ for the tangent spaces of $X$ at $A,B$, so that $\dim S_1(X)\le  2n+1-(h+x)$. Now for a third point $C\in H\cap X$, we see that $\dim <T_A\cup T_B>\cap T_C\ge h-1+x$ so that $\dim S_2(X)\le  3n+2-2(h+x)$  and so on: the conclusion $\dim S_k(X)\le n(k+1)+k-k(h+x)$ follows.\pas
For the last inequality, just call $h'$ the minimal value such that $\dim G_{h',k}(X)\le n(k+1)+(k-h')(h'+1)-x-2(h-h')$. If $h'=0$ then $\dim S_k(X)\le n(k+1)+k-x-2h$, otherwise we get
$\dim G_{h'-1,k}(X)\ge n(k+1)+(k-h'+1)h'-x-2(h-h')-1$ and the previous conclusion tells us that $\dim S_k(X)\le n(k+1)+k-k(h'+x+2(h-h'))\le n(k+1)+k-x-2h$.\qed
\enddemo

\remark{Example} For $h=1$, $k=2$ the previous Proposition yields: when $X$ is not a cone and $\dim S_1(X)=2n+1< r$, then $\dim G_{1,2}(X)<3n+2$ implies also $\dim S_2(X)<3n$. In this case, the statement is in fact also a consequence of Proposition 1.9 below. Compare with the Proposition 1.6, which just says $\dim S_2(X)\le 3n$.\endremark\smallskip

\proclaim{Theorem 1.8} Assume $S_{k-1}(X)\neq \pp r$ and $\dim G_{h,k}(X)<  (k-h)(h+1)+n(k+1)$. Then $\dim G_{h-1,k-1}(X)<(k-h)h+ nk$. In particular $\dim S_{k-h}(X)<n(k-h+1)+k-h$.
\endproclaim
\demo{Proof} Take $\pi\in G_k(X)$ general and consider a general $h$-plane $L\subset\pi$. By Corollary 1.3, we know that $\pi\cap X=\{p_0,\dots,p_k\}$. Let $L'$ be the intersection of $L$ with the span of $p_1,\dots,p_k$: it is a hyperplane in $L$ and it is also a general element of $G_{h-1,k-1}(X)$.\pas
Now move $\pi$ in the family of $(k+1)$-secant $k$-planes through $L$; the points $p_i$ move consequently, so also their spans move. If $L'$ moves with $\pi$, then it gives a family which is dense in $L$; thus any point of $L$ lies in some $k$-secant $(k-1)$-plane, which implies $S_k(X)=S_{k-1}(X)$, whence $S_{k-1}(X)=\pp r$ by Proposition 1.4, a contradiction. Thus $L'$ is fixed, hence it belongs to infinitely many elements of $G_{k-1}(X)$.\pas
The last statement follows observing that $S_{k-h}(X)\subset S_{k-h+1}(X)\subset\dots\subset S_{k-1}(X)$.
\qed\enddemo

We will apply the previous result to the case $k=2$, $h=1$. It reads:  $r>2n+1$ and $\dim G_{1,2}(X)<3n+2$ implies $\dim S_1(X)<2n+1$.

\proclaim{Theorem 1.9} Assume $\dim G_{1,2}(X)<3n+2$ and $\dim S_1(X)=2n<r$. Then $X$ is a cone.\endproclaim
\demo{Proof} Take a general 2-plane $\pi\in G_2(X)$; we may assume, by Corollary 1.3, $\pi\cap X=\{ A,B,C\}$. Take a general point $P$ in the line $<A,B>$ and fix a general line $L\subset\pi$,  through $P$.\pas
Any line in $\pi$ is an element of $G_{1,2}(X)$, so, by our assumptions, it is contained in an infinite family of 3-secant planes. Move $\pi $ in the family of 3-secant planes through $L$, we get a new plane $\pi'$, which cuts $X$ in $A',B',C'$; as observed in the previous proof, by continuity, the hypothesis $S_1(X)\neq\pp r$ implies that $P$ is still contained in the line $<A',B'>$; so $L$ induces in this way a 1-dimensional family of secant lines through $P$.\pas
Move now $L$ in $\pi$ and consider the induced families of secants. Since $P$ is general in $S_1(X)$, these families must coincide, for $P$  cannot be contained in a 2-dimensional family of secants, by our assumptions. It follows that for $L'\subset \pi$ general through $P$, there exists a 3-secant plane containing $L'$ and $<A',B'>$; but this implies that every plane through $A',B'$, in the 3-space $M$ spanned by $\pi,A',B'$, is 3-secant to $X$. In particular, $X$ meets $M$ in a curve $\Gamma'$, passing through $C$.\pas
Observe that $L$ is a general line of $\pi$ and $M$ is also the span of $\pi$ and the plane $\pi'$ of $L,A',B'$; also observe that $\pi'$ is 3-secant and the family of 3-secant planes through $L$ is 1-dimensional, by Proposition 1.6, since $\dim S_1(X)=2n$.
It follows that $M$ is the span of two general planes in the component of the family of 3-secant planes through $L$ containing $\pi$; this component is unique, since $\pi$ is general.\pas
If $\Gamma'$ is non degenerate, then by Proposition 1.5 it must be a cone, otherwise $S_1(\Gamma')$ fills up the whole of $M$ and $S_2(X)=S_1(X)\neq\pp r$, contradiction. Hence $\deg\Gamma'\le 2$. If $\Gamma'$ is a conic, then $X$ contains a conic through two general points. If $\Gamma$ is such a conic through $A,B$, then all the lines through $P$ in the plane of $\Gamma$ are secant to $X$. Since $P$ is general in $S_1(X)$, then among these lines there is also the span $<A',B'>$ for there is just one component of secants through $P$ containing $<A,B>$; on the one hand this is the family of the secants through $P$ in the plane of $\Gamma$, on the other hand it contains $<A',B'>$. Thus $M$ contains also a conic through $A,B$ and proposition 1.5 shows that $S_1(X)=S_2(X)$, a contradiction.\pas
We conclude that $\Gamma'$ is a line, i.e. $X$ contains a line through any general point.\pas
Change now $C$ and $B$ and change $P$ with the point $Q=L\cap<A,C>$. As above, we get that there exists a line in $X$ passing through $B$ and contained in the span of $\pi$ and $\pi'$, i.e. in $M$; hence $M$ contains three lines of $X$, through the points $A,B,C$; since these lines are not in $\pi$, they form a non-degenerate curve. Since $M$ is contained in $S_2(X)$, then by Proposition 1.5 the three lines form a cone, since $S_1(X)$ cannot contain $M$ by assumption and by Proposition 1.4.\pas
Hence we get that for any triple of points of $X$ there pass 3 lines meeting in a common point. As in the proof of Proposition 1.7, this is possible only if $X$ is a cone.
\qed\enddemo 

\heading 2. Surfaces and their Grassmannians $G_{12}$\endheading

Through this section, $S$ is always an integral surface in some projective space. We want to classify surfaces for which the variety $G_{1,2}(S)$ has dimension smaller than the expected value 8; in other words, a general line in $G_{1,2}(S)$ 
is contained in infinitely many 3-secant planes. Since the situation is clear in $\pp 4$, for reasons of dimension, we may assume that $S\subset\pp r$, $r\ge 5$. It turns out that there are few surfaces with this property; namely:

\proclaim{Theorem} The integral surfaces in $\pp r$, satisfying  $\dim G_{1,2}(S) <8 $ are either cones or  rational normal surfaces of (minimal) degree $4$ in $\pp 5$, but not Veronese surfaces.\endproclaim

\remark{Remark 2.1} It is well known that a general line in $\pp 5$ lies in some 3-secant plane to a Veronese surface $S$. Just to give an elementary argument, look at $S$ as the locus of conics of rank 1 and observe that any pencil of conics is contained in some net generated by conics of rank 1. It is classically known that a general projection of a Veronese surface to $\pp 3$ has 3 double lines, forming a cone.\endremark\smallskip

The proof of the theorem is divided in several steps. We shall use often the following classical Lemma, due to C.Segre:

\proclaim{Segre's Lemma} Let $S\subset\pp N$, $N\ge 4$ be a non degenerate integral surface containing a 2-dimensional family of plane curves. Then the curves have degree $\le 2$ and $S$ is either a Veronese surface or a projection of a Veronese surface in $\pp 4$.
\endproclaim
\demo{Proof} The original proof is in [S]. See [CS] or [M] for a modern proof.\qed\enddemo

\remark{Step 1} \it We may reduce ourselves to surfaces in $\pp 5$ with $S_1(S)=S_2(S)=\pp 5$.\pas
\rm Indeed it follows from Theorem 1.8 that surfaces in $\pp r$, $r>5$ for which $\dim G_{1,2}(S)<8$, also satisfy $\dim S_1(S)<5$. Since classically it is known that all surfaces in $\pp 6$ with this last property are cones, they are yet considered in the classification.\pas
Furthermore, by Theorem 1.9 all surfaces with $\dim G_{1,2}(S)<8$ and $\dim S_1(S)=4$ still are cones. Thus we may also assume that $S$ is not a Veronese surface.\qed
\endremark\smallskip

Consider now the incidence variety: 
$$I'\subset G_2(S)\times \pp 5\quad I'=\{(\pi,Q): Q\in\pi\}.$$
call $p,q$ the projections; by Proposition 1.1, $\dim I'=8$. Since $S_2(S)=\pp 5$, then $p$ dominates $\pp 5$; so for $P\in \pp 5$ general, all components of $p^{-1}(P)$ have dimension 3.\pas
Choose $P$ general and choose one component $L_P$ of $p^{-1}(P)$.

\remark{Step 2} \it Let $W_P=p(q^{-1}(q(L_P)))=$ the union of all planes belonging to $L_P$. Then $W_P$ is an irreducible variety containing $S$.\pas\rm
Indeed the irreducibility of $W_P$ follows immediately from the irreducibility of $L_P$.\pas
Assume $S\not\subset W_P$. Then the inverse image of $L_P$ in $S^3\dots\to G(2,5)$ dominates only a curve $\Gamma$. We get an irreducible component $Y$ of $\Gamma^3$ such that for $Q_1,Q_2,Q_3$ general in $Y$, their span lies in $L_P$, i.e. $P\in<Q_1,Q_2,Q_3>$ and, for dimensional reasons,  general elements of $L_P$ are obtained in this way.\pas
Write $Y=\Gamma_{01}\times\Gamma_{02}\times\Gamma_{03}$, with each $\Gamma_{0i}$ component of $\Gamma$ and take  $(A,B)$ general in $\Gamma_{01}\times\Gamma_{02}$. Since $P$ is general and $S_1(S)=\pp 5$, then $P$ lies only in finitely many secant lines to $S$, hence we may assume $P\notin<A,B>$; it follows that for $C\in\Gamma_{03}$ general, we have $<A,B,C>=<P,A,B>$ fixed. This contradicts the fact that $Y$ induces a 3-dimensional family of planes.
\qed\endremark\smallskip

\remark{Step 3}  $dim\  W_P = 4.$\pas
Assume $\dim W_P=5$. Then for $Q\in\pp 5$ general there exists $\pi\in L_P$ with $Q\in\pi$; then the line $<Q,P>$ belongs to $\pi$. Since the two points $P,Q$ are general, we see that a general line belongs to a 3-secant plane, contradicting the assumption $\dim G_{1,2}(S)<8$.\pas
Assume $\dim W_P\le 3$. Since $W_P$ is irreducible and contains both $S$ and the general point $P$, then $\dim W_P=3$. For $Q\in S$ general there exists a 3-secant plane contained in $W_P$ passing through $Q$; it follows that the line $<P,Q>$ lies in $W_P$, hence $W_P$ is the cone over $S$ with vertex $P$. Since $W_P$ contains a 3-dimensional family of planes, the projection of $W_P$ from $P$ is a surface containing a 3-dimensional family of lines: it is classical that such a surface cannot exist.\qed
\endremark\smallskip

Choose now two points $A,B\in S$ such that $P$ belongs to the line $\ell=<A,B>$; since $P$ is general, we may assume that $A,B$ are general in $S$. Define a rational map:
$$\Psi: S\dots\to G(2,5)$$
which sends $C\in S\setminus\{A,B\}$ to the plane $<A,B,C>$; call $L'_P$ the closure of the image. Clearly $L'_P$ is irreducible and by construction it lies in $q(p^{-1}(P))$.\smallskip\noindent

From now on, we start with $A,B$ general, then we construct $L'_P$ and choose the component $L_P$ of  $q(p^{-1}(P))$, containing $L'_P$.

\remark{Step 4} \it $dim\  L'_P=2$ and $p(q^{-1}(L'_P))=W_P$.\pas
\rm By Corollary 1.3, for $C\in S$ general we have $<A,B,C>\cap S=\{A,B,C\}$, hence $\Psi$ has finite general fibers, i.e. $\dim L'_P=2$.\pas
Now take $C\in S$ general; the fiber of $q^{-1}(L'_P)$ over $C$ is non empty (by construction) and finite, for otherwise $C$ lies in infinitely many planes of $L'_P$. Then $p(q^{-1}(L'_P))$ has dimension 4. Moreover for $Q$ general in $p(q^{-1}(L'_P))$, there exists $C\in S$ with $Q\in<A,B,C>\in L_P$; this means $p(q^{-1}(L'_P))\subset W_P$. The claim follows from irreducibility of $W_P$ and Step 3. 
\qed\endremark\smallskip

For $\pi\in L_P\setminus L'_P$, write $\Lambda_\pi$ for the span of $\pi$ and the line $\ell =<A,B>$.

\remark{Step 5} \it $\Lambda_\pi$ is a 3-dimensional space contained in $W_P$ and $W_P$ is the union of spaces $\Lambda_\pi$, as $\pi$ varies. Hence for $\pi\in L_P$ general, the intersection $\Lambda_\pi\cap S$ contains a curve. \pas\rm
Clearly $\dim\Lambda_\pi=3$ for $P\in\ell\cap\pi$. Also for $Q\in\pi$ general, by step 4, there is $C\in X$ with $Q\in<A,B,C>$, whence $<Q,A,B>\subset W_P$, so $\Lambda_\pi\subset W_P$. When $\pi$ varies, the corresponding $\Lambda_\pi$ define a variety of dimension at least 4, it coincides with $W_P$. Since $S\subset W_P$, for $\pi$ general it follows that $\Lambda_\pi\cap S$ has dimension 1.
\qed\endremark\smallskip

Observe that $\Lambda_\pi\cap S$ might have some isolated point, for the ambient fourfold $W_P$ might be singular somewhere. Thus although $A,B$ belong to $\Lambda_\pi\cap S$, unfortunately we are not allowed to conclude, a priori, that there exists a curve contained in $\Lambda_\pi\cap S$ and passing through $A,B$. This makes the argument more involved.\pas
Next step is crucial to override this difficulty and conclude the proof of the Theorem.

\remark{Step 6} \it Assume $S$ is not a cone. Then for $\pi,\pi'\in L_P$ general, we have $\Lambda_\pi\cap\Lambda_{\pi'}=\ell$.\pas
\rm Assume on the contrary that $\Lambda_\pi\cap\Lambda_{\pi'}$ is a plane $V$. First we show that this $V$ is fixed as $\pi,\pi'$ vary; indeed otherwise, by the Linear Lemma, the union of the spaces $\Lambda_\pi$ is a $\pp 4$ which contains $S$, contradicting the assumptions.\pas
The projection $S\dots\to \pp 2$ from $V$ is not dominant, for it contracts the curves $\Lambda_\pi\cap S$, which cover $S$. It follows that the projection of $S$ from the secant line $\ell\subset V$ is a cone in $\pp 3$. The conclusion now follows from:\endremark

\proclaim{Lemma} If the projection of an integral surface $S
\subset\pp N$ $N\ge 4$, from a general point $A\in S$, is a cone, then $S$ itself is a cone.\endproclaim
\demo{Proof} Call $\tau_A$ the projection; if $F$ is a ruling of the cone $\tau_A(S)$, then $\tau_A^{-1}(F)$ is a plane curve; thus $\tau_A$ determines a 1-dimensional family of plane curves.\pas
Assume these curves are not lines; then they determine their planes, so moving $A$ on $S$ also the family moves; it follows that $S$ is covered by a 2-dimensional family of plane curves; by Segre's Lemma, $S$ is a projection of a Veronese surface. On the other hand it is well known that projecting a Veronese surface or a cubic surface in $\pp 4$, which is not a cone, from a general point of the surface, we cannot get a cone.\pas
Assume that all curves $\tau_A^{-1}(F)$ are lines, so $S$ is ruled; since it is not a cone, its lines meet  the  line $R$ which joins $A$ to the vertex of $\tau_A(S)$, in a moving point; moving $A$ generically, we find a new line $R'$ which is intersected by all the lines of $S$; then $S$ lies in the span of $R,R'$, a contradiction.\qed\enddemo

Since cones are included in our classification, we may assume from now on that {\it $S$ is not a cone}.

\remark{Step 7} {\it $P$ is not contained in any tangent line to smooth points of $S$. Moreover $\ell$ is the unique secant line to $S$ contained in $\Lambda_\pi$ and passing through $P$.}\pas\rm
The first claim follows easily from the observation that the union of all tangent planes to regular points of $S$ is a 4-dimensional variety.\pas 
Assume that for $\pi\in L_P$ general there exists another secant line $\ell'$ with $P\in\ell'\subset\Lambda_\pi$.  $\ell'$ must be fixed as $\pi$ varies, since $P$ is a general point of $\pp 5$ and $S_1(S)=\pp 5$, so $P$ is contained in finitely many secant lines. But this contradicts the previous step, which tells that  the intersection of two general spaces $\Lambda_\pi\cap \Lambda_{\pi'}$ is $\ell$. \qed\endremark\smallskip

Call now $\Gamma_\pi$ the union of all 1-dimensional components of $\Lambda_\pi\cap S$. As observed in step 5, $\Gamma_\pi$ is non empty, but it may be different from the intersection $\Lambda_\pi\cap S$.

\remark{Step 8} {\it For $\pi$ general in $L_P$, there are no components of $\Gamma_\pi$ which are plane curves and contain $A$ or $B$.}\pas\rm
Assume there exists such a component $\gamma$ through $A$. Since $\ell$ is a general secant line, then it is not tangent to $S$, for $S_1(S)=\pp 5$. Thus $\ell$ cannot coincide with the embedded tangent space $T$ to $\gamma$ at $A\in S$. It follows that $\gamma$ moves, when $\pi$ varies in $L_P$, for otherwise a general intersection $\Lambda_\pi\cap \Lambda_{\pi'}$ contains $T$, contradicting step 6.\pas
Since $\gamma$ moves, then $A$ is contained in a 1-dimensional family of plane curves lying on $S$; since $A$ is general, then $S$ contains a 2-dimensional family of plane curves; this is impossible by  Segre's Lemma, since $S_1(S)=\pp 5$.\qed\endremark\smallskip

\remark{Step 9}\it The variety $W_P$ cannot contain a 2-dimensional family of spaces $\Lambda_P$. \pas\rm
Indeed $S$ is non degenerate and we have the following general:
\proclaim{Lemma} Let $Y$ be an irreducible variety of dimension $m\ge 2$, containing a 2-dimensional family $\Cal R$ of linear spaces of dimension $m-1$. Then $Y$ itself is a linear space.\endproclaim
\demo{Proof} The proof goes by induction on $m$, the case $m=2$ being classical (and an easy consequence of the Linear Lemma). For $m>2$, take a general hyperplane $H$ and consider the family $\{R_i\cap H: R_i\in\Cal R\}$ of subvarieties of $Y\cap H$. If this family is 2-dimensional, then we conclude by induction. Otherwise, for $R\in\Cal R$ general, the intersection $R_H=H\cap R$ (which is a general hyperplane in $R$) must be contained in infinitely many spaces of $\Cal R$; but in this case, varying $H$, the families $\{R'\in\Cal R:R_H\subset R'\}$ together dominate $\Cal R$, so $R$ meets another general element $R'\in\Cal R$ in a non-fixed subspace of codimension 1. It follows by the Linear Lemma that the total space of $\Cal R$ is a linear space of dimension $m$, contained in $Y$.\qed\enddemo

\remark{Step 10}\it We have $\deg\Gamma_\pi= 3$ and $\Gamma_\pi$ is a rational normal cubic passing through $A,B$.\pas\rm
Since $L_P$ is 3-dimensional and all planes of $L_P$ belong to some $\Lambda_\pi\subset W_P$, it follows by a dimensional count that the general plane of $\Lambda_\pi$ passing through $P$ is a 3-secant plane and by Corollary 1.3 it meets $S$ exactly in 3 points; this is clearly impossible unless $\Lambda_\pi\cap S$ is a curve of degree 3.\pas
Take now  a general plane $\pi'\subset\Lambda_\pi$ passing through $\ell$; it is a general 3-secant plane, so it meets $S$ in exactly 3 points, by Corollary 1.3. But $\deg\Gamma_\pi=3$ and $\Gamma_\pi\cap\ell\subset S\cap\ell=\{A,B\}$; also, by step 8, no components of $\Gamma_\pi$ through $A$ are plane curves. If $\Gamma_\pi$ does not contain (say) $A$, then a general plane through $\ell$ lying in $\Lambda_\pi$ meets $S$ in  more than 3 points. This is impossible by Corollary 1.3, for this plane is a general 3-secant plane to $S$. Thus $\Gamma_\pi$ is an irreducible smooth cubic passing through $A$ and $B$. \qed\endremark\smallskip 

\remark{Step 11}\it The linear system $|L|$ above determines a birational map from $S$  to a quadric surface in $\pp 3$. \pas\rm
Take $A,B\in S$ general; for all spaces $\Lambda_\pi$ as above, we find a rational normal cubic $\gamma \subset S$, passing through $A,B$ and contained in $\Lambda_\pi$; by step 6, moving $\pi$ we see that $S$ contains a 1-dimensional family of rational normal cubics which intersects transversally in $A,B$; moving $A,B$, we get a 3-dimensional family of rational normal cubics on $S$. Since these curves are contracted by the Albanese map, then $h^1\Cal O_S=0$ and they belong to a linear system $|L|$ of dimension 3; by Noether's theorem, $S$ is rational. Since by step 6 the intersection of two curves of $|L|$ through the general points $A,B$ is exactly $\{A,B\}$, then $L^2=2$ and the sequence:
$$0\to\Cal O_S\to\Cal O_S(L)\to\Cal O_\gamma(2)\to 0$$
implies $h^0\Cal O_S(L)=4$. We get that $|L|$ has no base points and it separates general points; the associated map thus sends $S$ birationally to a quadric in $\pp 3$.\qed\endremark\smallskip

Call $g:S\dots\to\pp 3$ the map associated to $|L|$.

\remark{Step 12: end of the classification} Assume $g(S)$ is a smooth quadric. The linear system of hyperplanes of $S\in \pp 5$ corresponds to some divisor class $(a,b)$ on $g(S)$;  since $g$ is biratonal on its image, $a,b>0$, moreover $a+b=3$: it follows that $S$  corresponds to the embedding of a smooth quadric surface via the linear system of type $(2,1)$ (or $(1,2)$).\pas
In fact, embedding a smooth quadric $Q$ with the linear system $|(2,1)|$ we find  a surface $S\subset \pp 5$ such that for a general points $A,B,C\in S$ there is a unique rational normal cubic $\gamma\subset S$ through them: $\gamma$ is the image of the unique conic through the corresponding points on $Q$. This curve $\gamma$ spans a $\pp 3$; call it $V$. If $R$ is any line contained in the plane spanned by $A,B,C$, then a general plane passing through $R$ and contained in $V$ is 3-secant to $S$. Hence $S$ is an example of a surface with few lines on 3-secant planes.\par
Assume $g(S)$ is a quadric cone. Call $X$ the blow up of $g(S)$ at the vertex and let $T$ be the class of the proper transform of a ruling; let $E$ be the exceptional divisor; the map $X\dots\to \pp 3$ corresponds to the class $2T+E$ on $X$. Call $aE+bT$ the class corresponding to the map $X\dots\to \pp 5$; then as above $(aE+bT)\cdot(2T+E)=3$ which implies $b=3$. Since $g$ is birational and $\dim(|3T|)=3$, we se that $a>0$; since $aT+bE$ has no fixed components and $E\cdot(aE+3T)=-2a+3$, we get $a=1$.\pas
Observe that $E+3T$ corresponds to  rational normal curves passing through the vertex of $g(S)$; $|E+3T|$ is very ample, in the blow up of $g(S)$ and, as above, this system determines a surface $S\subset\pp 5$ such that 3 general points of $S$ are contained in a rational normal cubic on it;  it follows that this embedding of quadric cones also provide example of surfaces in $\pp 5$ with few lines on 3-secant  planes. \qed\endremark\smallskip

\remark{Remark 2.2} \rm
Rational surfaces of degree 4 in $\pp 5$, except for the $2$-Veronese embedding of $\pp 2$, are  rational scrolls of degree 4. There are three types of such scrolls, corresponding to the rank two bundles $\Cal O(2)\oplus\Cal O(2)$, $\Cal O(1)\oplus\Cal O(3)$ and $\Cal O\oplus\Cal O(4)$ over the projective line. Scrolls of the first type can also be seen as a quadric of $\pp 3$ embedded in $\pp 5$ by means of the linear system $(2,1)$. The second type corresponds to the blowing up of a quadric cone at the vertex, embedded by the (pull-back of the) linear system of cubic curves through the vertex. All these scrolls $X$ are smooth and  have $\dim G_{1,2}(X)<8$, by step 12 of the theorem.\par
Scrolls of the third type are cones; also for such $X$ one has  $\dim G_{1,2}(X)<8$, for three general points $A,B,C$ span, together with the vertex, a 3-plane $L$ which meets $X$ in three lines, so every line $r$ in the 3-secant plane $<A,B,C>$ lies in infinitely many 3-secant planes: the planes in $L$ through $r$.\par
Observe that rational normal scrolls $X$ of degree 4 in $\pp 5$ are in fact classically known to project generically to $\pp 3$ as surfaces with no triple points; it follows that the general projection of $X$ to $\pp 4$ has no 3-secant lines through a general point, whence a general line of $\pp 5$ does not lie on any 3-secant plane to $X$.

\endremark
\smallskip 

\pas
Luca Chiantini - Dipartimento di Matematica - Via del Capitano, 15 - 53100 SIENA (Italy) - email: chiantini$\@$unisi.it\smallskip\pas
Marc Coppens - Department Industrieel Ingenieur en Biotechniek - Katholieke Hogeschool Kempen - Kleinhoefstraat 4 - B 2440 GEEL (Belgium) - email: marc.coppens$\@$khk.be

\end
\smallskip


\Refs\widestnumber\key{CS}

\ref\key B\by Bronowski J.\paper The sums of powers as simultaneous canonical espressions\jour Proc. Camb. Phyl. Soc.\yr 1933\pages 465-469\endref

\ref\key CS\by Ciliberto C., Sernesi E.\paper Singularities of the theta divisor and congruences of planes\jour J. Alg. Geom.\vol 1\yr 1992\pages 231-250\endref

\ref\key JH\by Harris J.\paper Algebraic Geometry, a first course\jour Springer  Graduate Text in Math.\vol 133\yr 1992\endref

\ref\key H\by Hartshorne R.\paper Algebraic Geometry\jour Springer Graduate Text in Math.\vol 52\yr 1977\endref

\ref\key M\by Mezzetti E.\paper Projective varieties with many degenerate subvarieties\yr 1994\jour Boll. UMI\vol 8B\pages 807-832\endref

\ref\key Sc\by Scorza G.\paper Determinazione delle varieta' a tre dimensioni di $S_r$, $(r\ge 7)$, i cui $S_3$ tangenti si tagliano a due a due\yr 1907\jour Rendiconti del Circolo Matematico di Palermo\vol 25\pages 193-204\endref  

\ref\key S\by Segre C. \paper Le superficie degli iperspazi con una doppia infinit\'a di curve piane o spaziali\jour Atti R. Accad. Sci. Torino\yr 1921\vol 57\pages 75-89\endref

\ref\key T\by Terracini A.\paper Sulla rappresentazione delle coppie di forme ternarie mediante somme di potenze di forme lineari\jour Ann. Mat. Pura Appl. \vol XXIV \yr 1915\pages 91-100 \endref

\ref\key Z\by Zak F.L.\paper Tangents and secants of algebraic varieties\jour Transl. Math. Monogr.\vol 127\yr 1993\endref

\endRefs
